\documentclass[11pt,a4paper]{article}
\usepackage{authblk}
\usepackage[utf8]{inputenc}
\usepackage[T1]{fontenc}
\usepackage{graphicx}
\usepackage{subfigure}
\usepackage{multirow}
\usepackage{amssymb}
\usepackage{amsmath}
\usepackage{amsfonts}
\usepackage[numbers,sort&compress]{natbib}
\usepackage{hyperref}
\usepackage{caption}
\usepackage{subfigure}
\usepackage{mathrsfs}
\usepackage[title]{appendix}
\usepackage{xcolor}
\usepackage{textcomp}
\usepackage{manyfoot}
\usepackage{booktabs}
\usepackage{algorithm}
\usepackage{algorithmicx}
\usepackage{algpseudocode}
\usepackage{listings}
\newtheorem{thm}{Theorem}[section]
\newtheorem{deff}{Definition}[section]

\newtheorem{ex}{Example}[section]
\numberwithin{deff}{section}
\numberwithin{thm}{section}
\title{Dynamical Analysis Of Fractional Order Generalized Logistic Map}
\author[1]{Sachin Bhalekar\footnote{Corresponding author email sachinbhalekar@uohyd.ac.in}}
\author[1]{Janardhan Chevala}
\author[2]{Prashant M. Gade}
\affil[1]{School of Mathematics and Statistics, University of Hyderabad, Hyderabad, 500046 India.}
\affil[2]{Department of Physics, Rashtrasant Tukadoji Maharaj Nagpur University, Nagpur, India.}
\date{}
\begin{document}
	
	\maketitle
	\begin{abstract}
		In this work, we propose a generalization to the classical logistic map. The generalized map preserves most properties of the classical map and has richer dynamics as it contains the fractional order and one more parameter. We propose the stability bounds for each equilibrium point. The detailed bifurcation analysis with respect to both parameters is presented using the bifurcation diagrams in one and two dimensions. The chaos in this system is controlled using delayed feedback. We provide some non-linear feedback controllers to synchronize the system. The multistability in the proposed system is also discussed.
	\end{abstract}
	
	\section{Introduction} 
	Generalizing classical mathematical science concepts is a critical term in research in this field. Such generalizations retain the properties of the original idea and provide flexibility to the models describing some physical phenomena. One of the most popular generalizations is the fractional calculus, which deals with the operators of generalized order \cite{podlubny1998fractional}. The order of fractional operators can be a non-integer, unlike those in the classical calculus. The fractional order derivatives have proved helpful in modeling memory and hereditary properties in various natural systems \cite{hilfer2000applications,magin2012fractional,mainardi2022fractional,kumar2023generalized,sarkar2023dynamic}.
	\par The maps are discrete-time dynamical systems arising in the real-life systems \cite{martelli2011introduction}. In \cite{may1976simple}, Robert May proposed the logistic map describing the population dynamics. The system gained popularity among the scientists due to its chaotic behavior.
	\par In their seminal work, Miller and Ross \cite{miller1989fractional} initiated the studies on the fractional difference calculus. Atici and Eloe \cite{atici2007transform,atici2009initial} generalized the sum and difference operators to include the real (non-integer) order. Abdeljawad \cite{abdeljawad2011riemann} defined the left and right Caputo fractional sums and differences.
	Wu and Baleanu \cite{wu2014discrete} presented the fractional order version of the classical logistic map and studied chaos and synchronization \cite{wu2014chaos}.
	\par In this article, we propose a generalized logistic map (GLM). The classical logistic map is given by $x(t+1)=g(x(t))$, where $g(z)=\mu z (1-z)$ \cite{may1976simple}. If $x(t+1)=g(x(t))$ is the given map, then we can write it as a difference equation $\Delta x(t)=g(x(t))-x(t)$. We propose a fractional order generalization as $\Delta^{\alpha} x(t)=f(x(t+\alpha-1))-x(t+\alpha-1),$ where $f(z)=g(z)/(1+rg(z))$, $r\in \mathbb{R}$ and $0<\alpha\leq 1$. This generalization is in two directions: order $\alpha$ and parameter $r$. We get the classical logistic map for order value $\alpha=1$ and the parameter value $r=0$. The complexity of this map is higher than that of its classical counterpart. It has memory properties due to fractional order and shows chaos for a wider range of parameter values. Zeros of $f$ are the same as those of $g$. The function $f$ remains bounded even if $g$ approaches to infinity.  The derivative of $f$ at any point $x$ is $f'(x)=\frac{g'(x)}{(1+r g(x))^2}$. Hence $f'(x)=0$ if and only if $g'(x)=0$. Thus, the extrema of the function $g$ are the same as those of function $f$. The map $f$ does not blow up, and iterates are bounded over the entire real line unless one starts at the pole $g(x)=\frac{-1}{r}$. Such rational maps defined over the entire real line have not been studied enough in the literature.
	\par The paper is organized as follows:
	In Section \ref{prel}, we provide a few key definitions and theorems. In Sections \ref{glm}, we discuss the stability analysis of each equilibrium point of a GLM. In Section \ref{global},  we sketch bifurcation diagrams for different values of fractional order $\alpha$ with various parameter values $\mu$, $r$, and describe the period-doubling case. In Section \ref{control}, we use delayed feedback to control chaos in GLM. In section \ref{syn.}, we design some controllers to synchronize a GLM and provide some examples. In Section \ref{mult.}, we discuss the multistability of GLM of integer and fractional order cases. Finally, the results are summarized in Section \ref{con.}.

	\section{Preliminaries} \label{prel}
	In this section, we present some basic definitions and results.\\ 
	Let $h>0$, \,$ a\in \mathbb{R}$,
	$(h\mathbb{N})_a=\{a,a+h,a+2h,\dots\}$ and $\mathbb{N}_a=\{a,a+1,a+2,\dots\}$.
	
	\begin{deff}(See\cite{bastos2011discrete,ferreira2011fractional,mozyrska2015transform})
		For a function $x:(h\mathbb{N})_a\longrightarrow \mathbb{R}$, the forward h-difference operator is defined as
		$$(\Delta_hx)(t)=\frac{x(t+h)-x(t)}{h},$$
		where $t\in(h\mathbb{N})_a$. \\
		Throughout this paper, we take $a=0$ and $h=1$. 
	\end{deff}
	
	\begin{deff}\cite{mozyrska2015transform}
		For a function $x:\mathbb{N}_0\longrightarrow \mathbb{R}$, the fractional sum of order $\beta >0$\, is given by 
		$$(\Delta^{-\beta} x)(t)=\frac{1}{\Gamma(\beta)} \sum_{s=0}^n{\frac{\Gamma(\beta+n-s)}{\Gamma(n-s+1)}}x(s),$$  
		where $t=\beta +n$, $n\in \mathbb{N}_0$.
	\end{deff}
	
	\begin{deff}\cite{fulai2011existence,mozyrska2015transform}  
		Let $\mu>0$ and $m-1<\mu<m$, where $m\in\mathbb{N}.$
		The $\mu$th fractional Caputo-like difference is defined as
		$$\Delta^\mu x(t)=\Delta^{-(m-\mu)}(\Delta^mx(t)),$$
		where $t\in\mathbb{N}_{m-\mu}$ and
		$$\Delta^mx(t)=\sum_{k=0}^m \left(\begin{array}{c}m\\k\end{array}\right)(-1)^{m-k}x(t+k).$$
	\end{deff}
	
	\begin{thm}\cite{fulai2011existence}
		The difference equation $$\Delta^{\alpha} x(t)=f(x(t+\alpha-1))-x(t+\alpha-1),$$ where $0<\alpha \leq 1$, $t \in \mathbb{N}_{1-\alpha}$, 
		is  equivalent to
		\begin{equation} 
			x(t)=x(0)+\frac{1}{\Gamma(\alpha)} \sum_{j=1}^t \frac{\Gamma(t-j+\alpha)}{\Gamma(t-j+1 )} (f(x(j-1))-x(j-1)), \label{1}
		\end{equation}
		where $t \in \mathbb{N}_{0}$.
	\end{thm}
	
	\begin{deff}\cite{elaydi2005systems,hirsch2012differential}
		A steady state solution  $x_*$ of (\ref{1}) is a real number satisfying $f(x_*) = x_*$.
	\end{deff}

	\begin{deff} \cite{elaydi2005systems,hirsch2012differential}
		An equilibrium $x_*$ is stable if for each  $\epsilon>0$, there exists $\delta>0$ such that $|x_0 -x_* | < \delta $ implies
		$|x(t) - x_* | < \epsilon$, $t=1,2,3,...$\\
		If $x_*$ is not stable then it is unstable.
	\end{deff}
	
	\begin{deff} \cite{elaydi2005systems,hirsch2012differential}
		An equilibrium point $x_*$  is asymptotically stable  if it is stable and there exists $\delta>0$ such that $|x_0 -x_* | < \delta $ implies $ lim_{t\to\infty}x(t) =x_*$.
	\end{deff}
	
	\subsection{Stability analysis} \label{stab.}
	We linearize  Eq.(\ref{1}) near an equilibrium $x_*$ with the help of first-order Taylor's approximation. Therefore, we get 
	\begin{equation} 
		y(t)=y(0)+\frac{1}{\Gamma(\alpha)} \sum_{j=1}^t \frac{\Gamma(t-j+\alpha)}{\Gamma(t-j+1 )} ((a-1) y(j-1)), \label{2}
	\end{equation}
	where $a=f'(x^*)$ and $y(t)=x(t)-x^*$ is infinitesimally perturbed solution. Thus the local stability properties of the fractional order non-linear map Eq.(\ref{1}) in a neighborhood of an equilibrium $x_*$ are the same as those of the linear map Eq.(\ref{2}).
	\begin{thm} \label{thm1}
		\cite{buslowicz2013necessary,stanislawski2013stability,vcermak2015explicit,pakhare2022synchronization}
		An equilibrium point $x_*$ of the non-linear map (\ref{1}) is asymptotically stable if $1-2^\alpha <f'(x^*) <1$.
	\end{thm}
	
	\section{A GLM and its stability} \label{glm}
	Consider the generalized logistic map (GLM) defined by (\ref{1}) with
	$$f(x)=\frac{\mu x(1-x)}{1+r \mu x(1-x)},$$
	where $\mu$ and $r$ are control parameters. Note that, this system reduces to the classical logistic map\cite{may1976simple} $x(t)=f(x(t-1))$ if $r=0$ and $\alpha=1$.\\
	In this paper, we study the stability analysis of a fractional order map (\ref{1}).
	
\begin{figure}[H]
		\centering
		\includegraphics[width=0.8\textwidth]{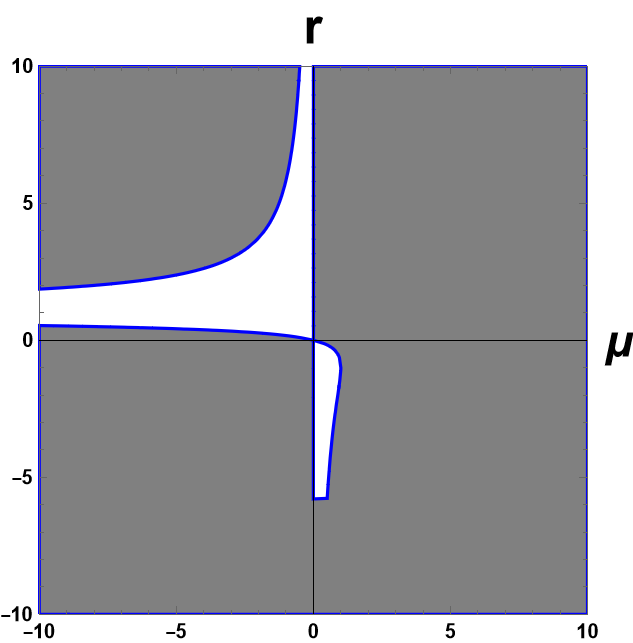}
		\caption{ The region of existence for $x_2^*$ and $x_3^*$.}
		\label{fig1}
	\end{figure}

	We have 
	$f(x)= \frac{\mu x(1-x)}{1+r \mu x(1-x)}$ and $f'(x)=\frac{\mu - 2 x \mu}{(-1 + r (-1 + x) x \mu)^2}$ . The fixed points of f give equilibrium points of (\ref{1}). We get $x_1^*=0$, $x_2^*=\frac{\mu + r \mu - 
		\sqrt{\mu} \sqrt{4 r + \mu - 2 r \mu+ r^2 \mu}}{2 r \mu}$ and  $x_3^*=\frac{\mu + r \mu +
		\sqrt{\mu} \sqrt{4 r + \mu - 2 r \mu+ r^2 \mu}}{2 r \mu}$ as equilibrium points. The equilibrium point $x_1^*$ is independent and other equilibrium points  $x_2^*$, and $x_3^*$ depend on parameters $\mu$ and $r$. Note that, the equilibrium points $x_2^*$ and $x_3^*$ are real when $4 r \mu + \mu^2 - 2 r \mu^2 + r^2 \mu^2 \ge 0$. The region of existence of $x_2^*$ and $x_3^*$ is given in the Fig. (\ref{fig1}).
	
	\subsection{Stability of $x_1^*$. }
	We have $f'(x_1^*)=\mu$. The stability of $x_1^*$ is independent of the parameter r. Therefore, the stable region in $\mu r$- plane is a vertical strip $1-2^{\alpha} < \mu <1 $ where $0< \alpha <1$. 
	
	\subsection{Stability of $x_2^*$. }
\begin{figure}[H]
		\centering
		\subfigure[Stable region of $x_2^*$ for $\alpha=0.2$.]
		{\includegraphics[height=2.0in,keepaspectratio,width=2.0in]{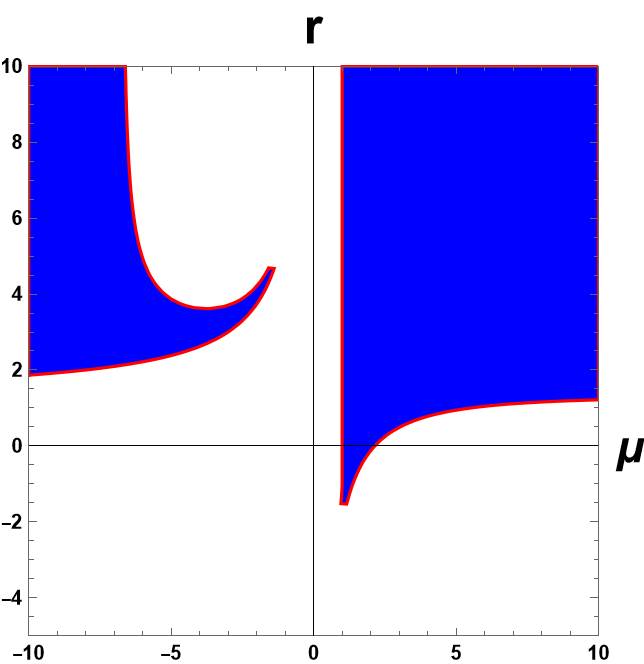}}\hspace{0.2cm}
		\subfigure[Stable region of $x_2^*$ for $\alpha=0.5$.]
		{\includegraphics[height=2.0in,keepaspectratio,width=2.0in]{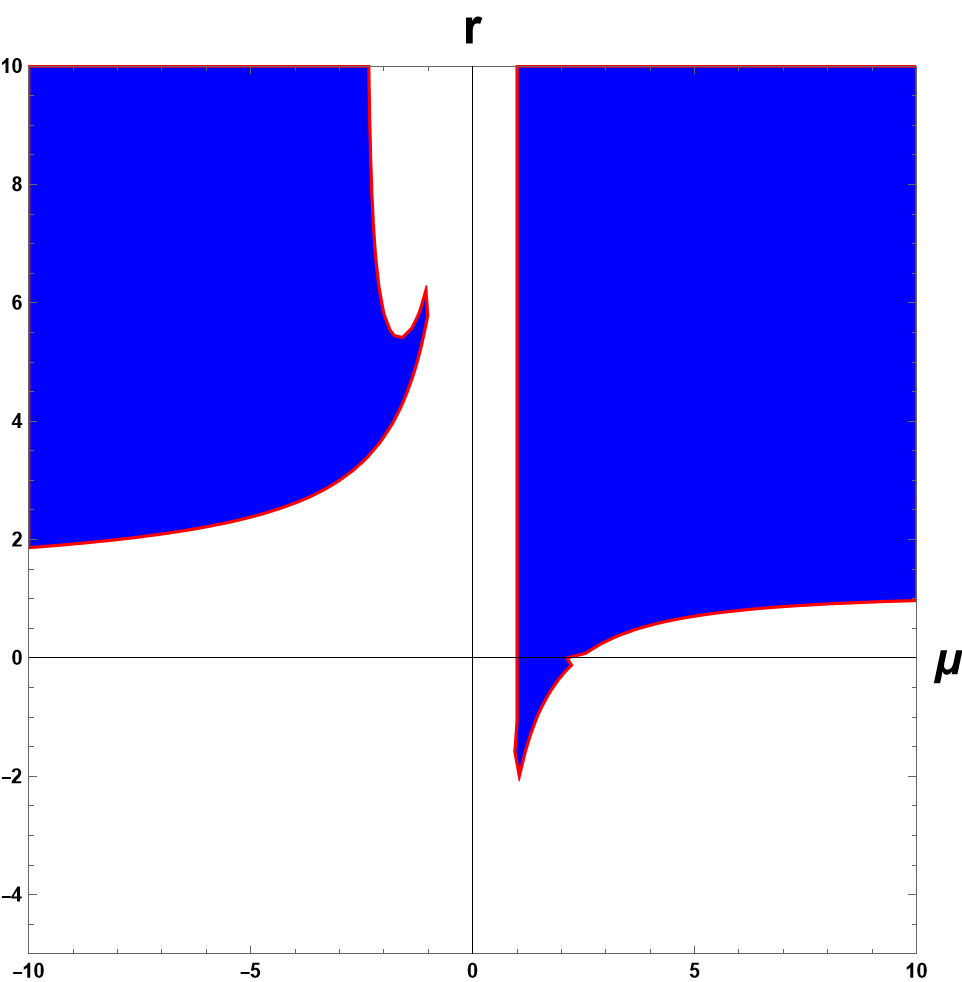}}\hspace{0.3cm}
		\subfigure[Stable region of $x_2^*$ for $\alpha=0.8$.]
		{\includegraphics[height=2.0in,keepaspectratio,width=2.0in]{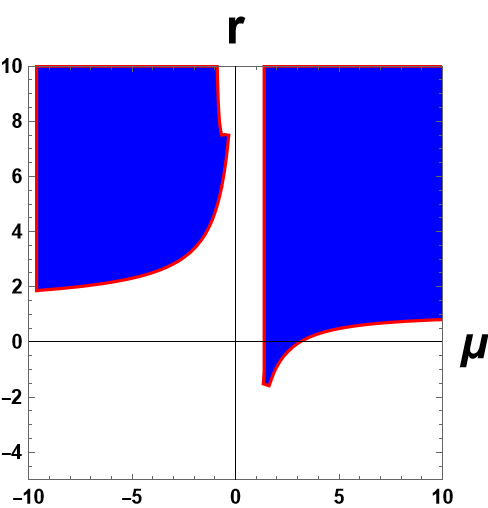}} 
		\caption{The stability regions of $x_2^*$ with different values of $\alpha$.} \label{fig2}
	\end{figure}

	We have $f'(x_2^*)= \frac{4 r(- \sqrt{\mu} + \sqrt{4 r + \mu - 2 r \mu + r^2 \mu )}}{\sqrt{\mu} ((-1 + r) \sqrt{\mu} + \sqrt{
			4 r + \mu - 2 r \mu + r^2 \mu})^2}$. \\
	The condition for stability of $x_2^*$ is $1-2^{\alpha} <  f'(x_2^*) <1 $ where $0< \alpha <1$.\\
	The region of stability of  $x_2^*$ in $\mu r-$ plane for $\alpha=0.2$,  $\alpha=0.5$ and $\alpha=0.8$ is shown in Fig. (\ref{fig2}).

	\subsection{Stability of $x_3^*$. }
\begin{figure}[H]
		\centering
		\subfigure[Stable region of $x_3^*$ for $\alpha=0.2$.]{\includegraphics[height=2.0in,keepaspectratio,width=2.0in]{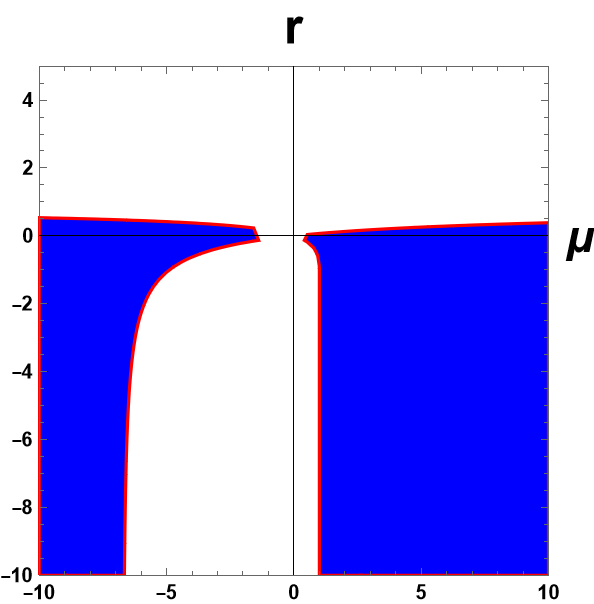}}\hspace{0.2cm}
		\subfigure[Stable region of $x_3^*$ for $\alpha=0.5$.]{\includegraphics[height=2.0in,keepaspectratio,width=2.0in]{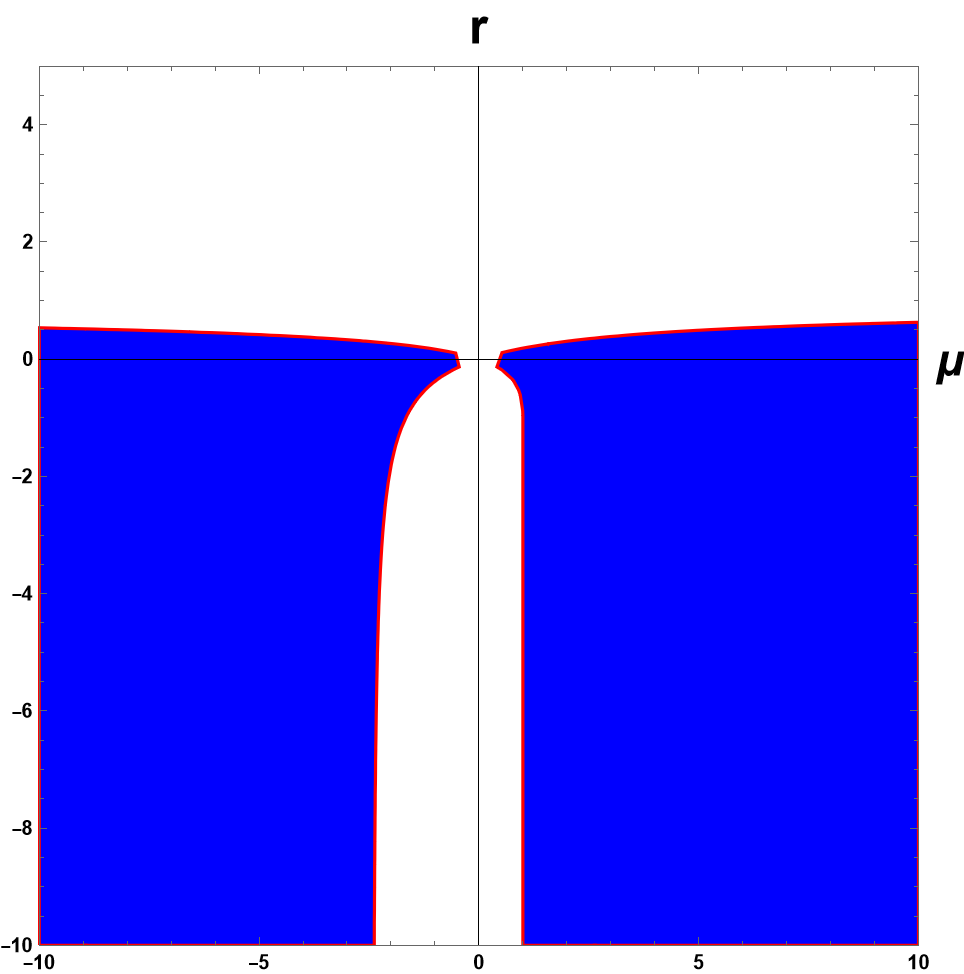}}\hspace{0.3cm}
		\subfigure[Stable region of $x_3^*$ for $\alpha=0.8$.]{\includegraphics[height=2.0in,keepaspectratio,width=2.0in]{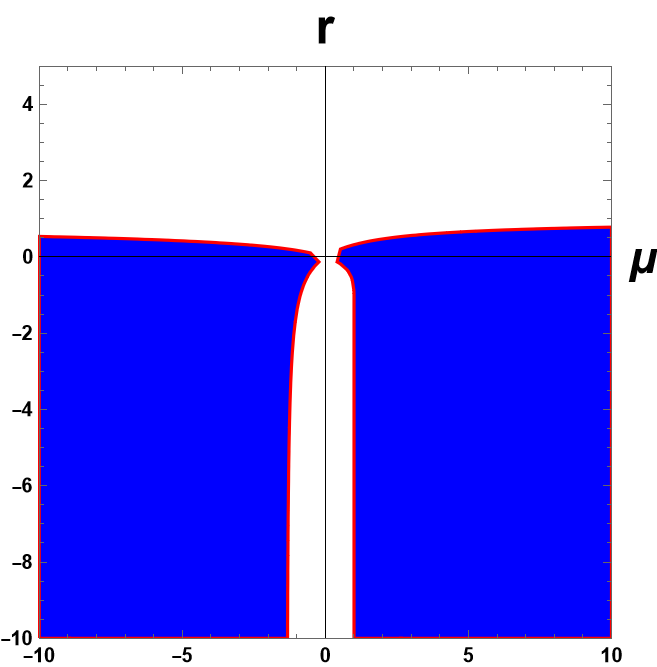}} 
		\caption{The stability regions of $x_3^*$ with different values of $\alpha$.} \label{fig3}
	\end{figure}

	We have $f'(x_3^*)= -\frac{4 r( \sqrt{\mu} + \sqrt{4 r + \mu - 2 r \mu + r^2 \mu)}}{\sqrt{\mu} ((-(-1 + r) \sqrt{\mu}) + \sqrt{
			4 r + \mu - 2 r \mu + r^2 \mu})^2}$. \\
	The condition for stability of $x_2^*$ is $1-2^{\alpha} <  f'(x_3^*) <1 $ where $0< \alpha <1$.\\
	The region of stability of  $x_2^*$ in $\mu r-$ plane for $\alpha=0.2$, $\alpha=0.5$ and $\alpha=0.8$ is shown in Fig. (\ref{fig3}).

	\section{Global bifurcation and period-doubling} \label{global}
\begin{figure}[H]
		\centering
		\subfigure[Global bifurcation diagram for $\alpha=0.2$.]{\includegraphics[height=4.0in,keepaspectratio,width=4.0in]{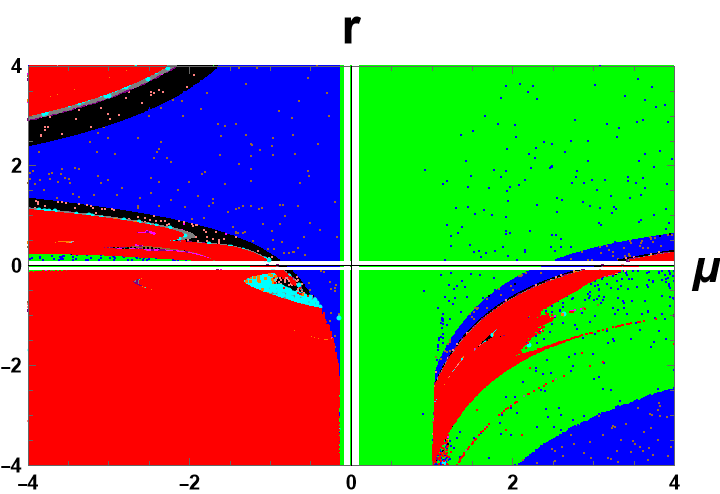}}\hspace{0.5cm}
		\subfigure[Global bifurcation diagram for $\alpha=0.8$.]{\includegraphics[height=4.0in,keepaspectratio,width=4.0in]{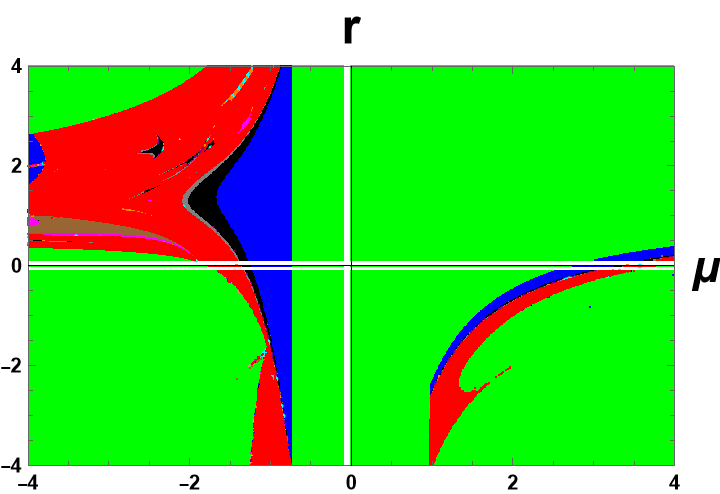}} 
		\caption{The global bifurcation diagrams with different values of $\alpha$.} \label{fig4}
	\end{figure}

	In Figs. \ref{fig4}(a) and \ref{fig4}(b), we sketch the phase diagrams for $\alpha=0.2$ and $\alpha=0.8$ respectively, in $\mu r$-plane showing different behaviors of the system (\ref{1}) with the color codes: green (asymptotic period-1), blue (asymptotic period-2), black (asymptotic period-4), gray (asymptotic period-8), and red (chaos). We observed from the numerical results in Fig. (\ref{fig4}) that, the system (\ref{1}) undergoes period-doubling bifurcation and shows chaos in each quadrant in $\mu r$-plane. Note that, we do not have exact periodic orbits for the fractional order systems \cite{bhalekar2023fractional}. We observed that the stability region increases in size with an increasing  $\alpha$ due to its dependence on $\alpha$.  
	
\begin{figure}[H]
		\centering
		\includegraphics[width=0.8\textwidth]{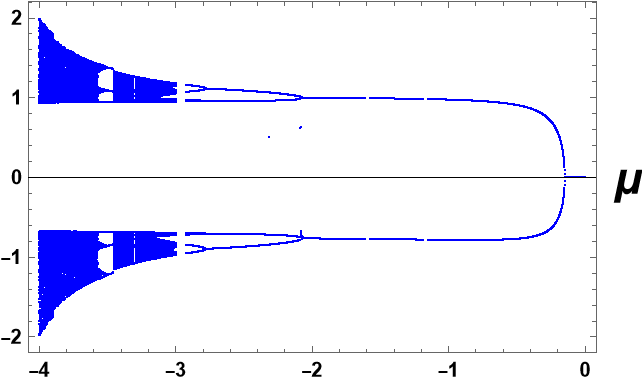}
		\caption{Bifurcation diagram at $r=3.5$, and $\alpha=0.2$. }
		\label{fig5}
	\end{figure}

	If we fix $r=3.5$, $\alpha=0.2$ and sketch the bifurcation diagram ( cf. Fig. (\ref{fig5}) ) for $\mu \in [-4,0]$ then we can observe that $\mu \in (-0.15,0)$ gives asymptotic period-1. We get asymptotic period-2 for  $\mu \in (-2.07,-0.15)$,  asymptotic period-4 for  $\mu \in (-2.77,-2.07)$, and asymptotic period-8 for  $\mu \in (-2.93,-2.77)$. Chaos is observed for $\mu \in (-4,-3.57)\cup(-3.45,-3.31)\cup(-3.29,-3.12)$. Periodic windows are observed for $\mu \in (-3.57,-3.45) \cup (-3.31,-3.29)$. 
	
	In this case, the Fig. \ref{fig6}(a) shows the orbit converging to $x_1^*$ for $\mu=-0.05$. For $\mu=-0.94$, we show an asymptotic period-2 orbit in Fig. \ref{fig6}(b). The asymptotic period-4 and asymptotic period-8 orbits for $\mu=-2.24$ and $\mu=-2.84$ respectively are given in Fig. \ref{fig6}(c) and Fig. \ref{fig6}(d). Fig. \ref{fig6}(e) shows a chaotic trajectory for $\mu=-3.85$.
	
	\begin{figure}[H]
		\centering
		\subfigure[Asymptotic period-1, $\mu=-0.05$.]{\includegraphics[height=2.0in,keepaspectratio,width=2.0in]{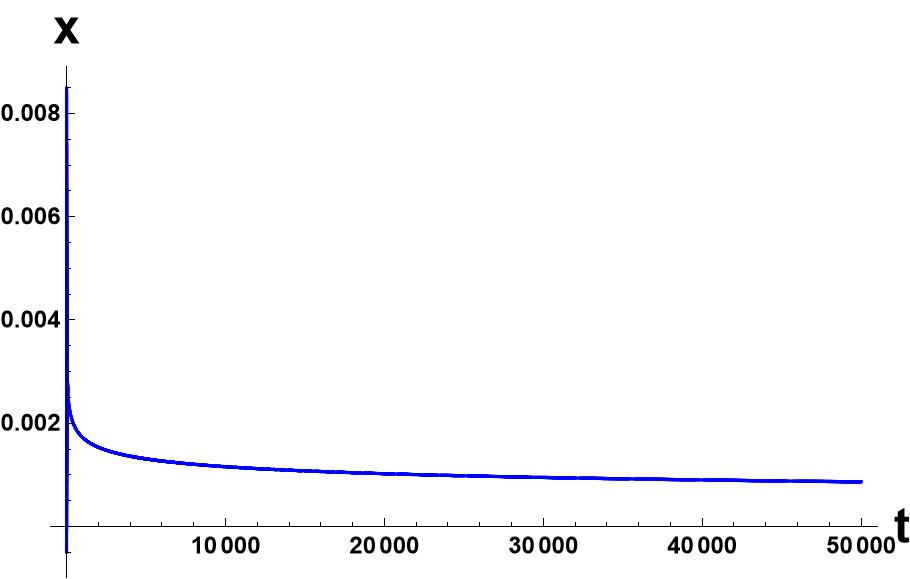}}\hspace{0.2cm}
		\subfigure[Asymptotic period-2, $\mu=-0.94$.]{\includegraphics[height=2.0in,keepaspectratio,width=2.0in]{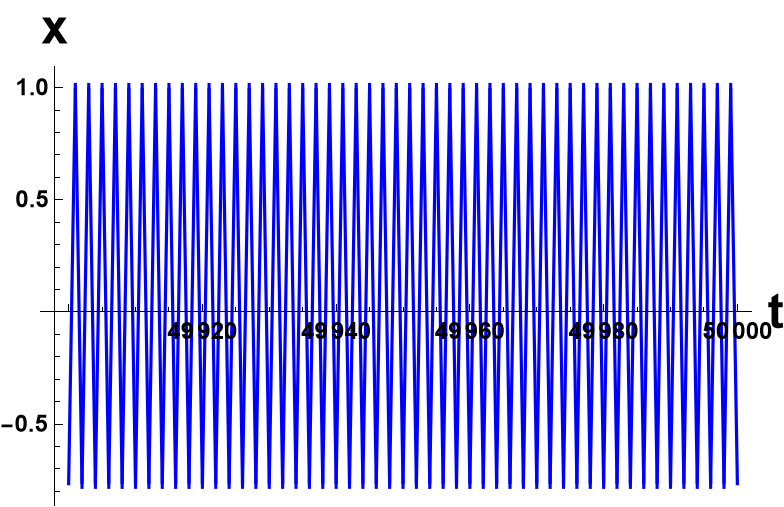}}\hspace{0.3cm}
		\subfigure[Asymptotic period-4, $\mu=-2.24$.]{\includegraphics[height=2.0in,keepaspectratio,width=2.0in]{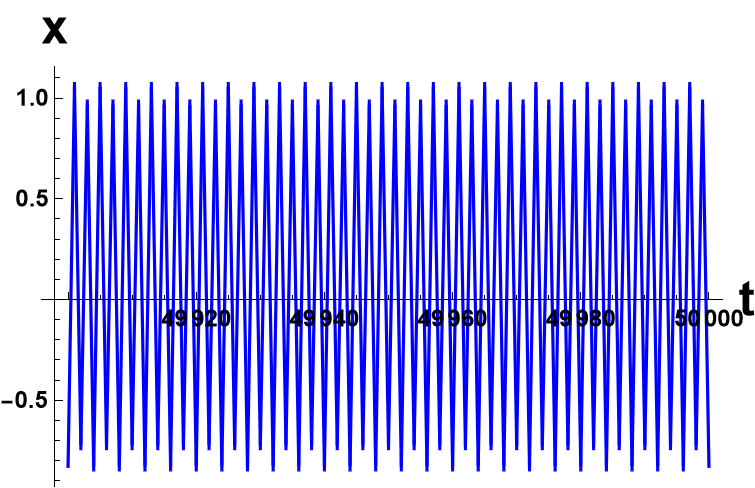}}
	\end{figure}
	
	\begin{figure}[H]
		\centering
		\subfigure[Asymptotic period-8, $\mu=-2.84$.]{\includegraphics[height=2.0in,keepaspectratio,width=2.0in]{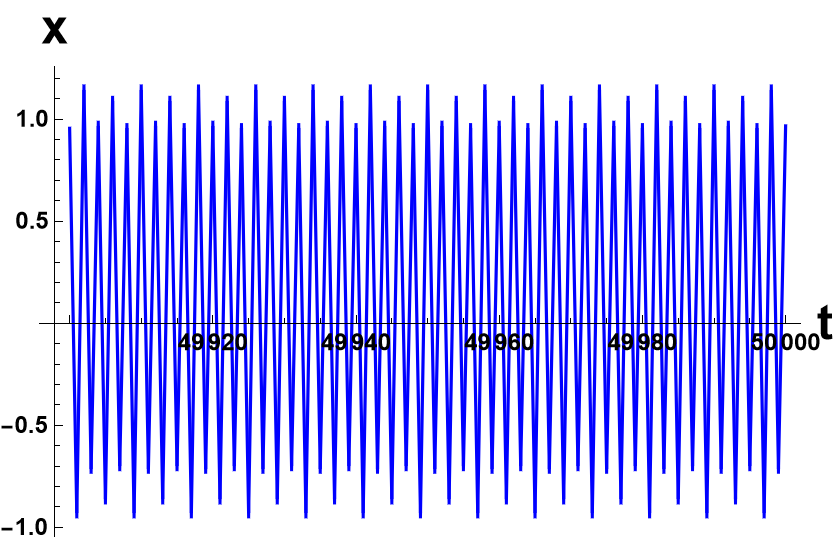}}\hspace{0.3cm}
		\subfigure[Chaos, $\mu=-3.85$.]{\includegraphics[height=2.0in,keepaspectratio,width=2.0in]{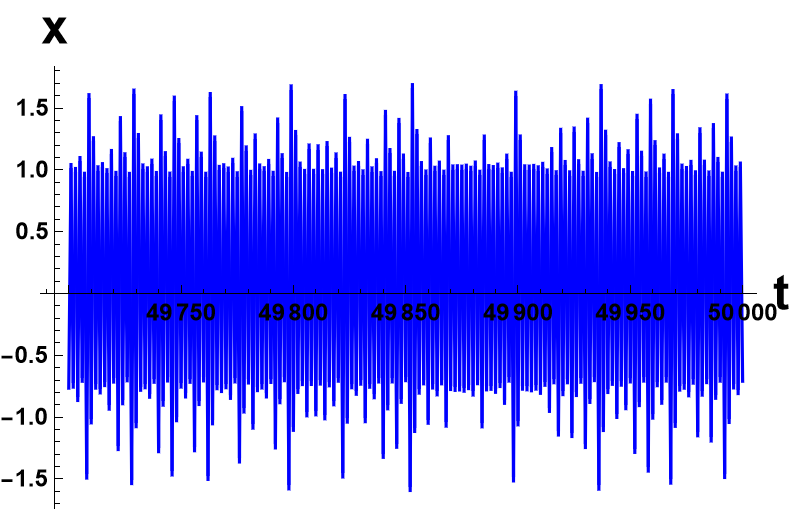}}
		\caption{Different behaviors of (\ref{1}) with different values of $\mu$ for fixed $r=3.5$ and $\alpha=0.2$.} \label{fig6}
	\end{figure}

	
	We also sketch the bifurcation diagrams for fixed $\alpha=0.2$ and various values of $r$ viz. $r=0.2$, $r=1$, $r=-0.3$ and $r=-2$ in Figs. \ref{fig7}(a)-\ref{fig7}(d) respectively. Here, the dynamics have the transitions from asymptotic stable point to asymptotic period-2, to asymptotic period-4, to asymptotic period-8, to asymptotic multi-period, then to chaos. 
	
	\begin{figure}[H]
		\centering
		\subfigure[$r=0.2$.]{\includegraphics[height=2.0in,keepaspectratio,width=2.0in]{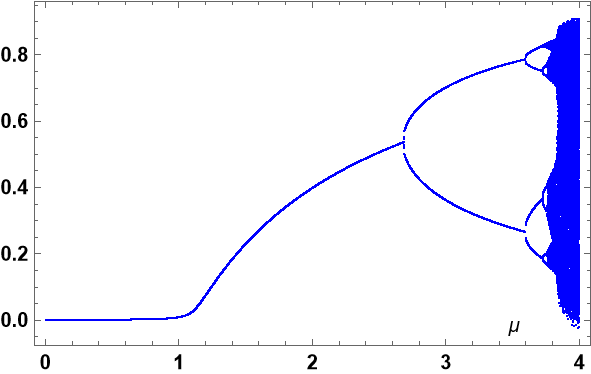}}\hspace{0.3cm}
		\subfigure[$r=1$.]{\includegraphics[height=2.0in,keepaspectratio,width=2.0in]{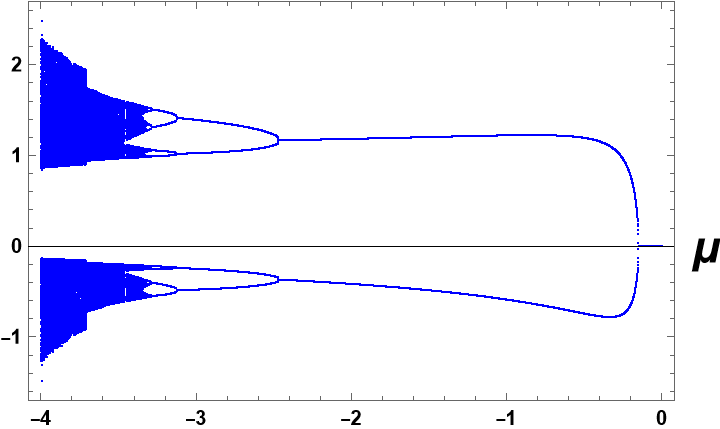}}
	\end{figure}
	
	\begin{figure}[H]
		\centering
		\subfigure[$r=-0.3$.]{\includegraphics[height=2.0in,keepaspectratio,width=2.0in]{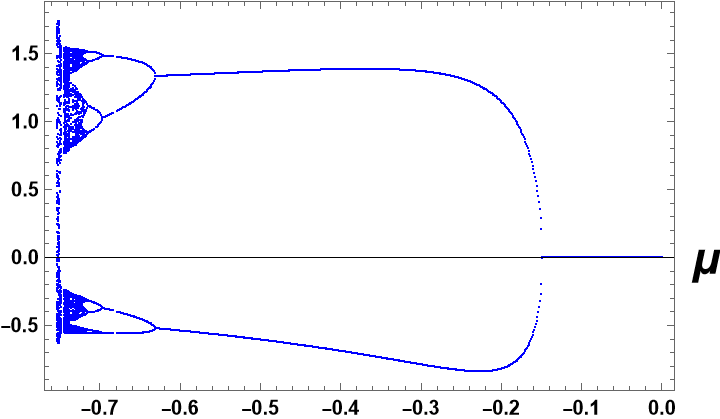}}\hspace{0.3cm}
		\subfigure[$r=-2$.]{\includegraphics[height=2.0in,keepaspectratio,width=2.0in]{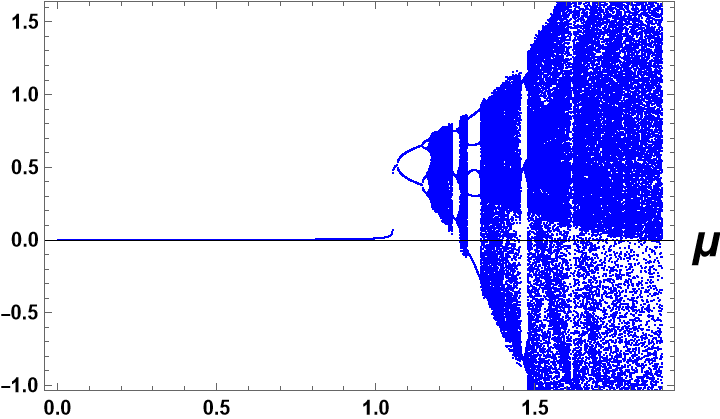}}
		\caption{Bifurcation diagrams with different values of $r$ for fixed $\alpha=0.2$.} \label{fig7}
	\end{figure}

\begin{figure}[H]
		\centering
		\includegraphics[width=0.8\textwidth]{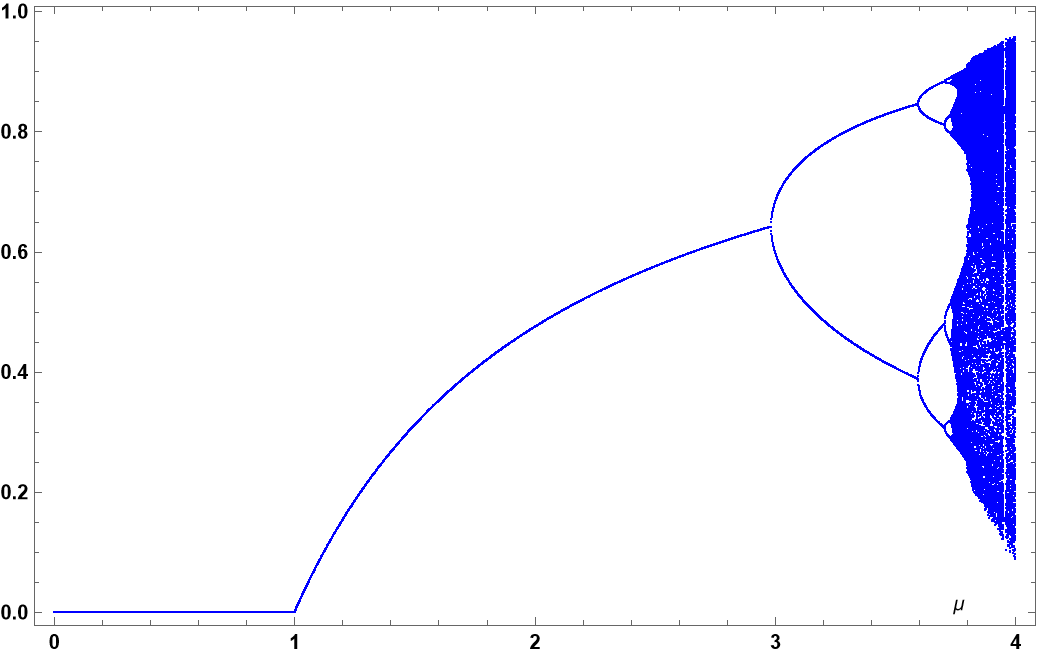}
		\caption{Bifurcation diagram at $r=0.1$, and $\alpha=0.8$. }
		\label{fig8}
	\end{figure}

	In the same way, we fix $r=0.1$, $\alpha=0.8$ and sketch the bifurcation diagram ( cf. Fig. (\ref{fig8}) ) for $\mu \in [0,4]$ to see that $\mu \in (0,3)$ provides asymptotic period-1. For  $\mu \in (3,3.58)$,  asymptotic period-2 is obtained; for  $\mu \in (3.58,3.71)$, asymptotic period-4 is obtained; and  for  $\mu \in (3.71,3.73)$, asymptotic period-8 is obtained. Chaos is seen for $\mu \in (3.76,4)$.

	In this instance, the orbit converges to $x_1^*$ for $\mu=2$ in Fig. \ref{fig9}(a). The figures \ref{fig9}(b), \ref{fig9}(c) and \ref{fig9}(d) display the asymptotic period-2, Period-4 and  period-8 orbits for $\mu = 3.4$, $\mu=3.6$ and $\mu=3.72$ respectively. For $\mu=3.9$, Fig. \ref{fig9}(e) depicts a chaotic trajectory.

	\begin{figure}[H]
		\centering
		\subfigure[Asymptotic period-1, $\mu=2$.]{\includegraphics[height=2.0in,keepaspectratio,width=2.0in]{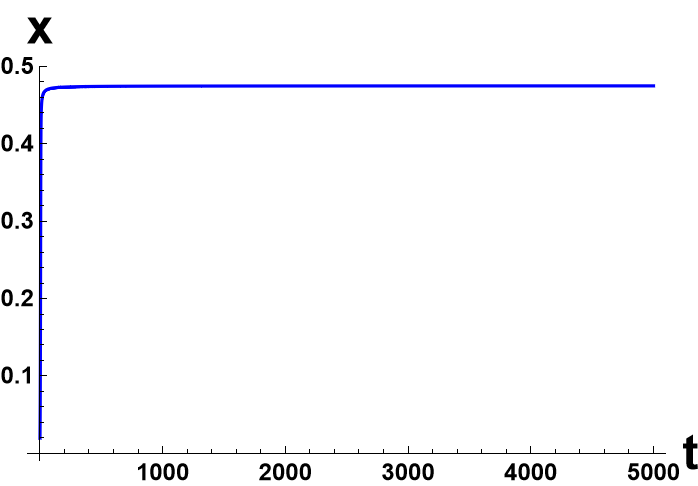}}\hspace{0.2cm}
		\subfigure[Asymptotic period-2, $\mu=3.4$.]{\includegraphics[height=2.0in,keepaspectratio,width=2.0in]{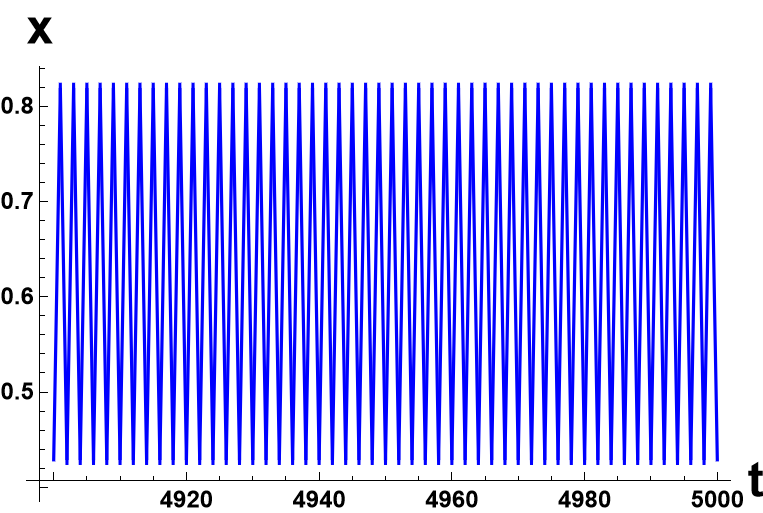}}\hspace{0.3cm}
		\subfigure[Asymptotic period-4, $\mu=3.6$.]{\includegraphics[height=2.0in,keepaspectratio,width=2.0in]{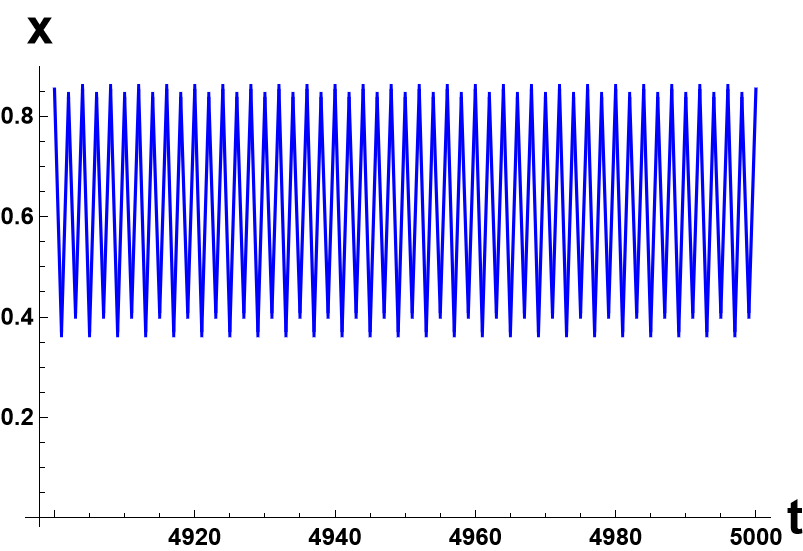}}
	\end{figure}
	
	\begin{figure}[H]
		\centering
		\subfigure[Asymptotic period-8, $\mu=3.72$.]{\includegraphics[height=2.0in,keepaspectratio,width=2.0in]{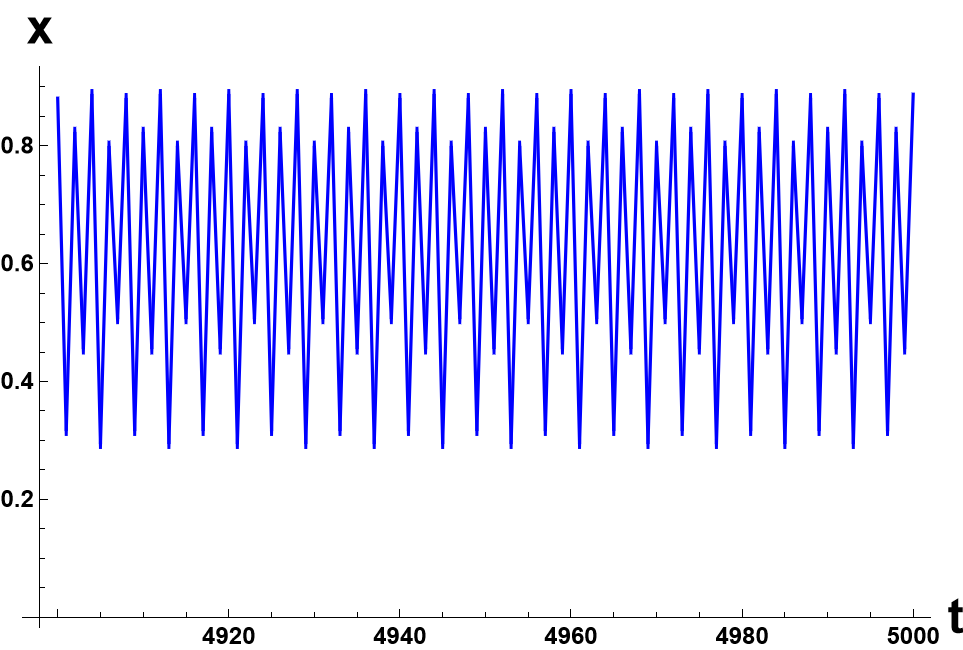}}\hspace{0.3cm}
		\subfigure[Chaos, $\mu=3.9$.]{\includegraphics[height=2.0in,keepaspectratio,width=2.0in]{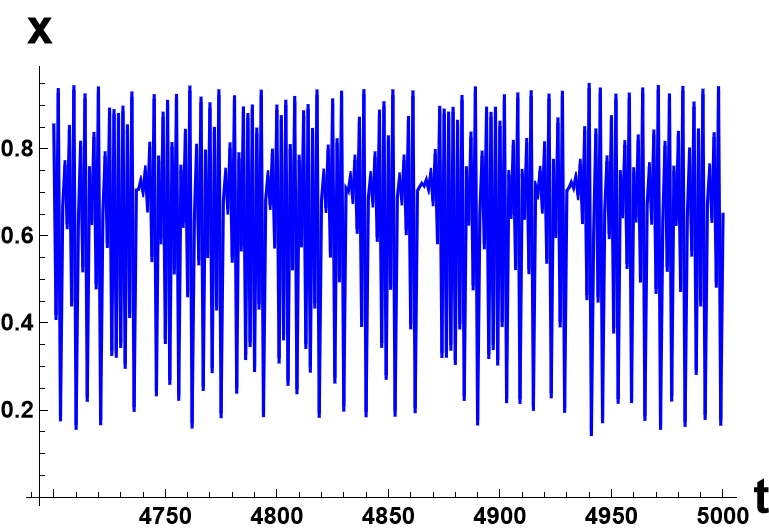}}
		\caption{Different behaviors with different values of $\mu$ for fixed $r=0.1$ and $\alpha=0.8$.} \label{fig9}
	\end{figure}
	
	Additionally, we draw the bifurcation diagrams fixed $\alpha=0.8$ for $r=0.1$, $r=1$, $r=-0.5$ and $r=-0.4$ in Figs. \ref{fig10}(a)-\ref{fig10}(d) respectively.
	\begin{figure}[H]
		\centering
		\subfigure[$r=0.1$.]{\includegraphics[height=2.0in,keepaspectratio,width=2.0in]{1Q1Bif_08_0.1.png}}\hspace{0.3cm}
		\subfigure[$r=1$.]{\includegraphics[height=2.0in,keepaspectratio,width=2.0in]{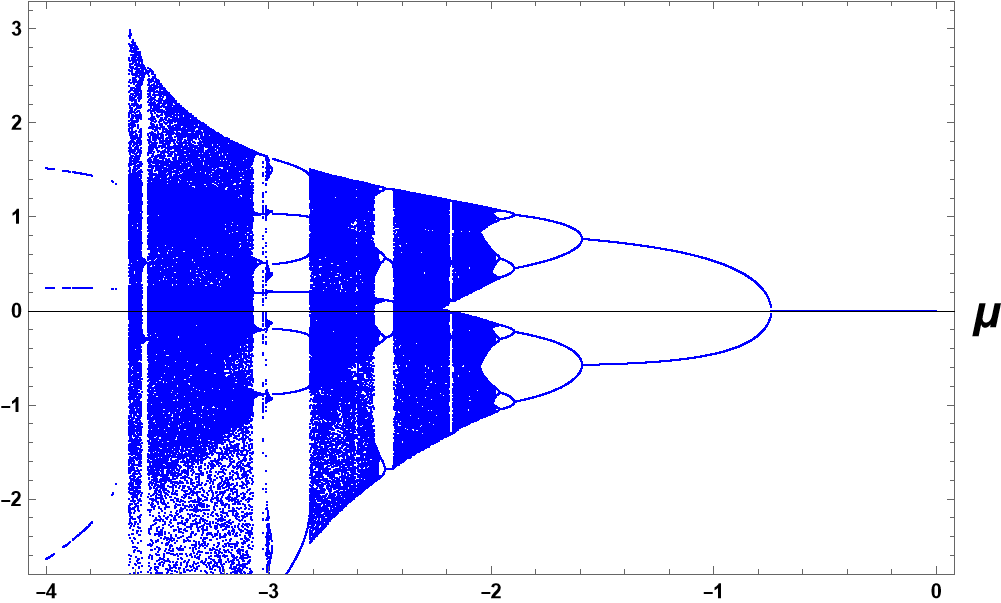}}
	\end{figure}
	
	\begin{figure}[H]
		\centering
		\subfigure[$r=-0.5$.]{\includegraphics[height=2.0in,keepaspectratio,width=2.0in]{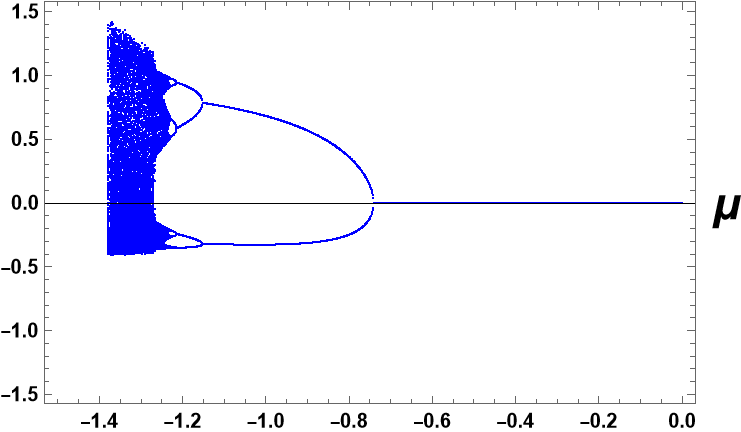}}\hspace{0.3cm}
		\subfigure[$r=-0.4$.]{\includegraphics[height=2.0in,keepaspectratio,width=2.0in]{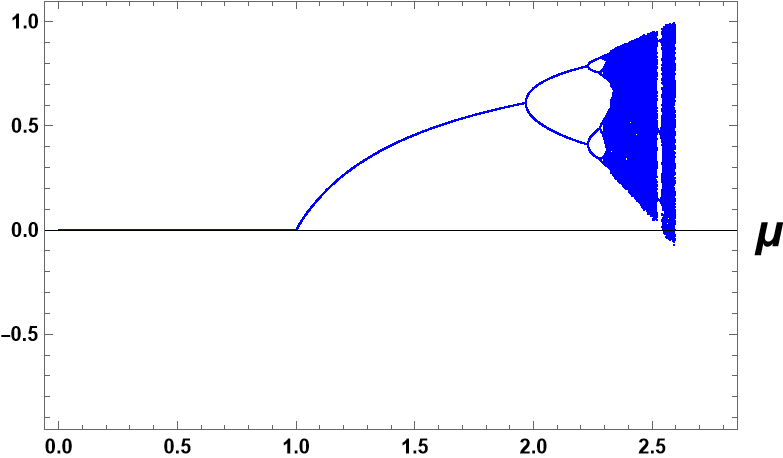}}
		\caption{Bifurcation diagrams with different values of $r$ for fixed $\alpha=0.8$.} \label{fig10}
	\end{figure}
	
	
	
	
	Moreover, some bifurcation diagrams are drawn for fixed $\alpha=0.5$ in Figs. \ref{fig11}(a)-\ref{fig11}(d) when $r=0.1$, $r=0.5$, $r=-0.3$ and $r=-1$ in that order.

	\begin{figure}[H]
		\centering
		\subfigure[$r=0.1$.]{\includegraphics[height=2.0in,keepaspectratio,width=2.0in]{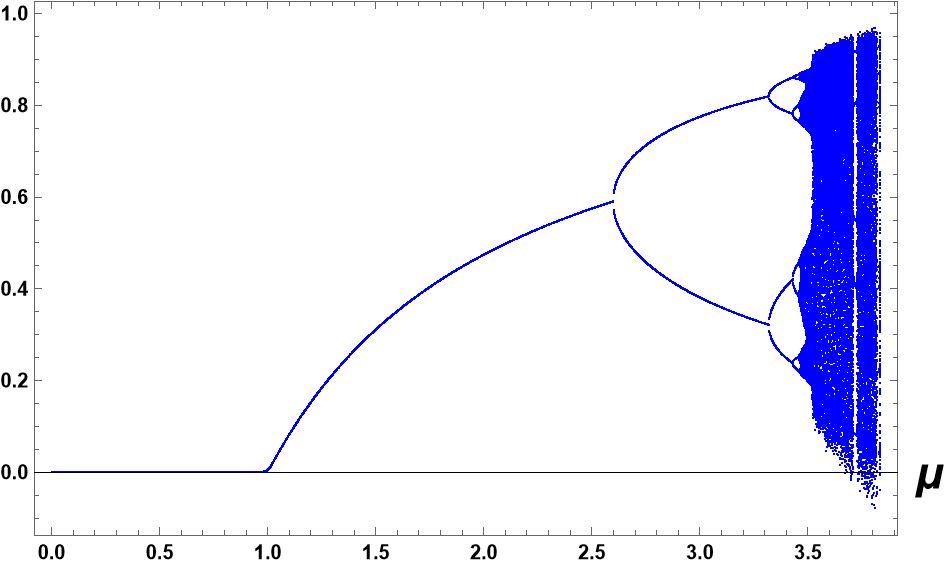}}\hspace{0.3cm}
		\subfigure[$r=0.5$.]{\includegraphics[height=2.0in,keepaspectratio,width=2.0in]{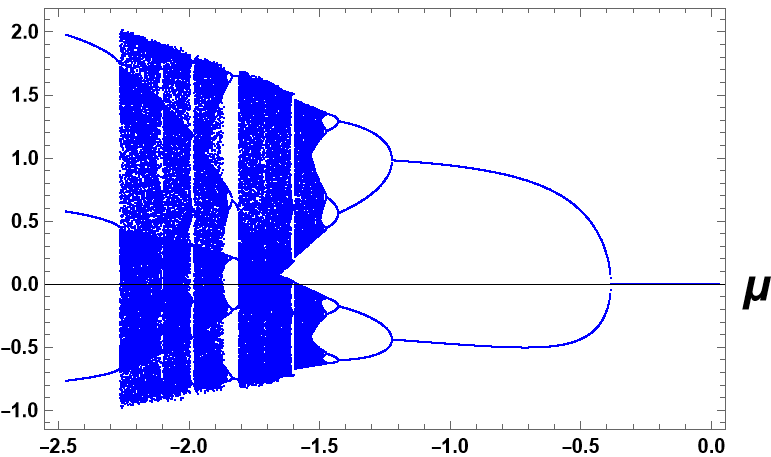}}
	\end{figure}
	
	\begin{figure}[H]
		\centering
		\subfigure[$r=-0.3$.]{\includegraphics[height=2.0in,keepaspectratio,width=2.0in]{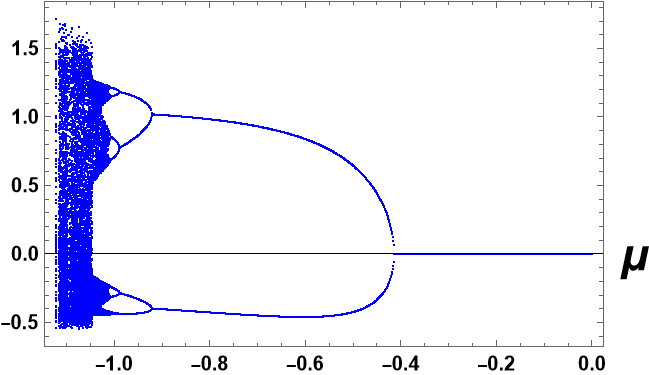}}\hspace{0.3cm}
		\subfigure[$r=-1$.]{\includegraphics[height=2.0in,keepaspectratio,width=2.0in]{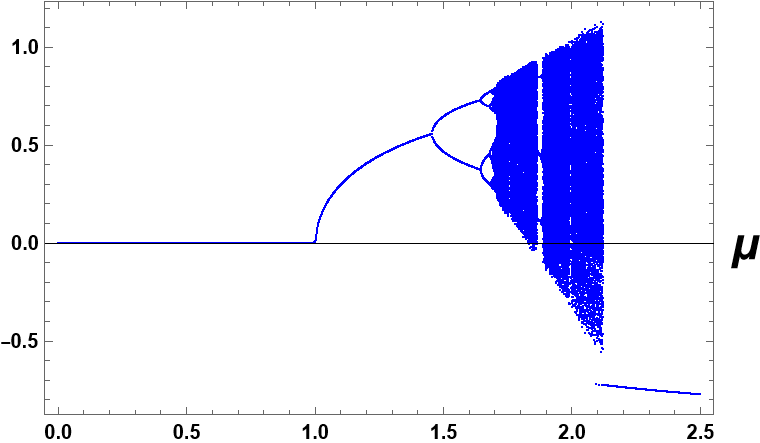}}
		\caption{Bifurcation diagrams with different values of $r$ for fixed $\alpha=0.5$.} \label{fig11}
	\end{figure}
	
	\section{Control of chaos using delayed feedback} \label{control}
	In this section, we propose a chaos-control in the system (\ref{1}) by using the delayed feedback control method derived in \cite{joshi2023controlling}. The controlled system is defined as
	\begin{equation}
		x(t)=x(0)+b x(t-\tau)+\frac{1}{\Gamma(\alpha)} \sum_{j=1}^t \frac{\Gamma(t-j+\alpha)}{\Gamma(t-j+1 )} (f(x(j-1))-x(j-1)), \label{3} 
	\end{equation}
	where $b\in \mathbb{R}$ is a control parameter and $\tau$ is the delay. We set $\tau=1$.
	In a neighborhood of the equilibrium point $x_1^*=0$, we have $f(x) \approx f(x_1^*)+(x-x_1^*)f'(x_1^*)=ax$, where $a=f'(0)=\mu$. The stable region of equilibrium point $x_1^*=0$ of the system (\ref{3}) is bounded by the curves \cite{joshi2023controlling} $a=1, a=1-2^{\alpha}(1+b)$ and the parametric curve 
	$$a(t)=2^{\alpha}(\sin{(\frac{t}{2})})^{\alpha} (\cos{(\frac{\alpha \pi}{2} + t(1-\frac{\alpha}{2}))}-\sin{(\frac{\alpha \pi}{2} + t(1-\frac{\alpha}{2})}) \cot{(\frac{\alpha \pi}{2}-\frac{t \alpha}{2})}+1,$$
	
\begin{figure}[H]
		\centering
		\subfigure[Stable region of fixed points at $\alpha=0.2$.]{\includegraphics[height=2.0in,keepaspectratio,width=2.0in]{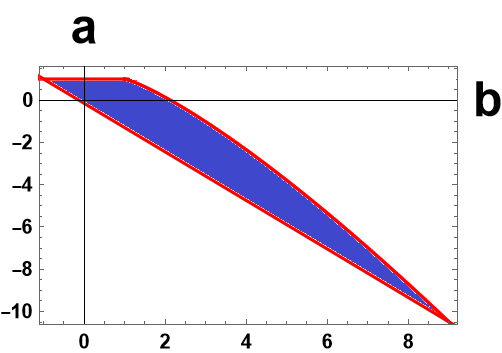}}\hspace{0.2cm}
		\subfigure[Stable region of fixed points at $\alpha=0.5$.]{\includegraphics[height=2.0in,keepaspectratio,width=2.0in]{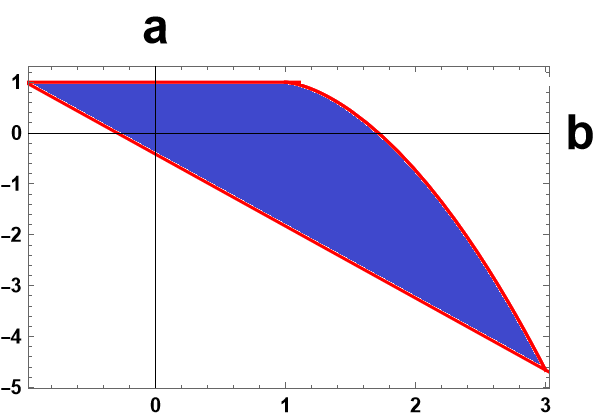}}\hspace{0.3cm}
		\subfigure[Stable region of fixed points at $\alpha=0.8$.]{\includegraphics[height=2.0in,keepaspectratio,width=2.0in]{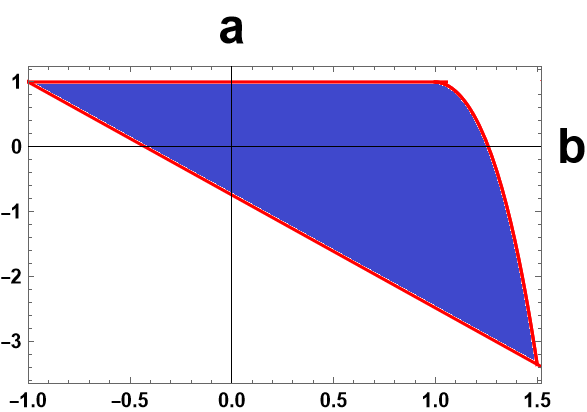}} 
		\caption{Stable fixed points of fractional GLM lie within the region enclosed by the $b-a$ curve with different $\alpha$ values .} \label{fig:12}
	\end{figure}

	$b(t)=\frac{\sin{(\frac{\alpha \pi}{2} + t(1-\frac{\alpha}{2}))}}{\sin{(\frac{\alpha \pi}{2}-\frac{t\alpha}{2})}},$ $t\in [0,2 \pi]$ in the b-a plane.  Thus the system (\ref{3}) can be controlled if the point $(b,a)$ with $a=f'(0)=\mu$ lies inside the region bounded by the above curves.\\
	For $\tau=1$, we have plotted this stable region in Figs. \ref{fig:12}(a)-\ref{fig:12}(c) for $\alpha=0.2, \alpha=0.5$ and $\alpha=0.8$ respectively.  
	
\begin{figure}[H]
		\centering
		\subfigure[b=0, $\alpha=0.2$, $r=3.5$.]{\includegraphics[height=2.0in,keepaspectratio,width=2.0in]{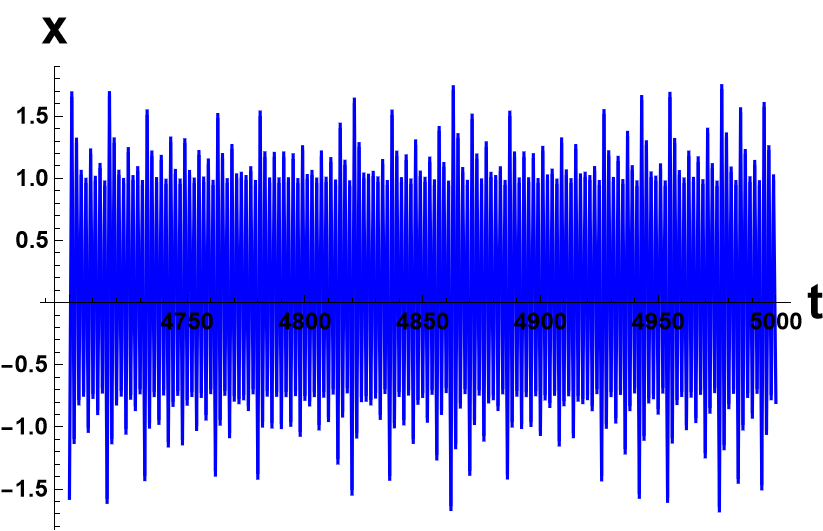}}\hspace{0.3cm}
		\subfigure[b=3.3, $\alpha=0.2$, $r=3.5$.]{\includegraphics[height=2.0in,keepaspectratio,width=2.0in]{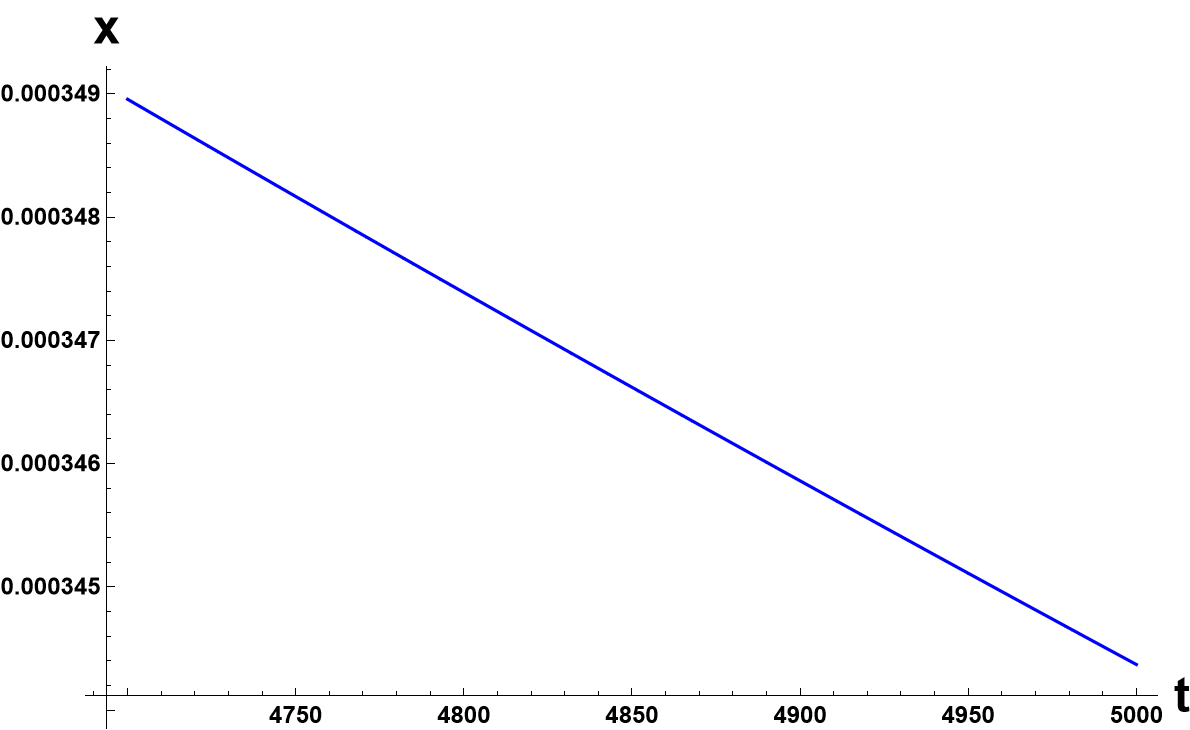}}
		
		\subfigure[b=0, $\alpha=0.8$, $r=2$.]{\includegraphics[height=2.0in,keepaspectratio,width=2.0in]{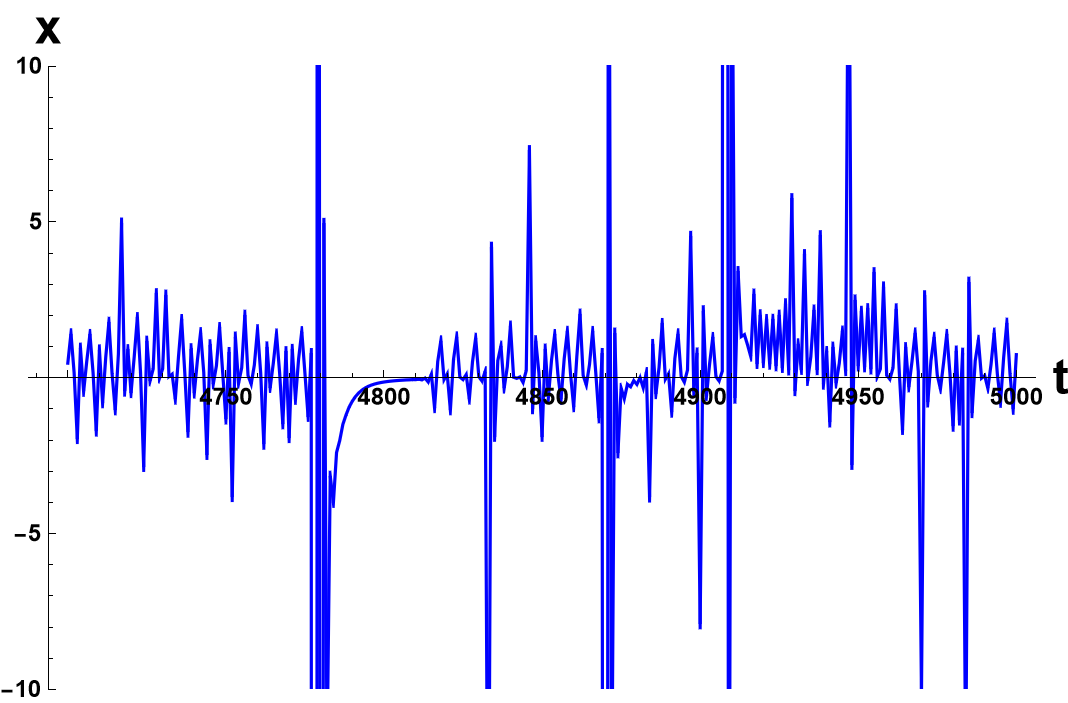}}\hspace{0.3cm}
		\subfigure[b=1.4, $\alpha=0.8$, $r=2$.]{\includegraphics[height=2.0in,keepaspectratio,width=2.0in]{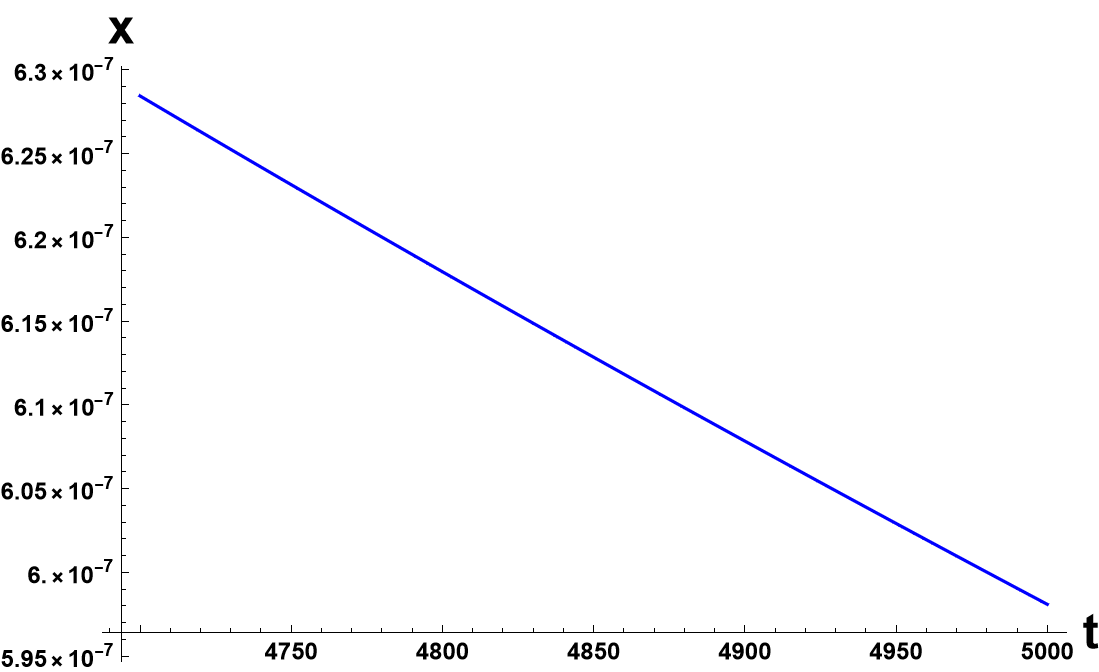}}
		\caption{Chaos control in system (\ref{1}) using delayed feedback.} \label{fig:13}
	\end{figure}

	If the system (\ref{1}) is chaotic for some set of parameters $\alpha, \mu$ and $r$ then we can control the chaos by selecting appropriate values of $(b,a)$ in the stable region. e.g. if $a=\mu=-3.85$, then $b\in (3.22,5.03)$ controls the chaos in the system. The chaotic solution for $\alpha=0.2$ and $b=0$ ( cf. Fig. \ref{fig:13}(a) ) are controlled by setting $b=3.3$ ( cf. Fig. \ref{fig:13}(b) ) in a stable region.
	
	Similarly, we control the chaos observed for $\alpha=0.8,\mu=-3$, and $r=2$ ( cf. Fig. \ref{fig:13}(c) ) by setting $b=1.4$ ( cf. Fig. \ref{fig:13}(d) ), where chaos control range is $b \in (1.29,1.48)$.
	
	\section{Synchronization of chaos} \label{syn.}
	The trajectories of two systems can oscillate together if we apply some control. This phenomenon is called synchronization \cite{pecora1997fundamentals,pecora2015synchronization}, which means going together. Since chaotic systems depend sensitively on the initial data, achieving synchronization in such systems may be challenging. In their pioneering work, Pecora and Carroll \cite{pecora1990synchronization} successfully developed a method showing that chaotic systems can also be synchronized. There are various methods viz. linear feedback control \cite{yassen2005controlling}, non-linear feedback control \cite{fu2012adaptive}, sliding-mode control \cite{kuz2022forced}, adaptive control \cite{saeed2023new} etc., to synchronize these systems. Bhalekar and Daftardar-Gejji \cite{bhalekar2010synchronization} developed an active control method to synchronize the fractional-order chaotic systems. Various other works include \cite{li2003synchronization,peng2007synchronization,lin2011h,yuan2012parameter,laarem2021new,qi2022synchronization}. However, there is insufficient literature on the synchronization of fractional order maps \cite{wu2014chaos,gade2021fractional,pakhare2022synchronization}.
	This section uses the non-linear feedback control and proposes the master-slave type synchronization in the fractional order GLM.

	\par We consider system (\ref{1}) as master system and the controlled slave system as
	\begin{equation}
		y(t)=y(0)+\frac{1}{\Gamma(\alpha)} \sum_{j=1}^t \frac{\Gamma(t-j+\alpha)}{\Gamma(t-j+1 )} \left[f(y(j-1))-y(j-1)+H(x(j-1),y(j-1))\right],\label{4}
	\end{equation}
	where H is the control function for the chaos synchronization.
	
	\begin{thm} \label{thm2}
		The master system (\ref{1}) gets synchronized with the slave system (\ref{4}) if  we use any one of the following controllers :\\
		(a) $H_{1}(x,y)= \frac{p(x-y)-p(x^2-y^2)}{(1+r \mu x (1-x))(1+r \mu y (1-y))}+k(x-y)$,\\
		(b) $H_{2}(x,y)= \frac{p(x-y)(1-2x)}{(1+r \mu x (1-x))(1+r \mu y (1-y))}+k(x-y)$,\\
		under the conditions $p=\mu$ and $-1<k<-1+2^{\alpha}$. 
	\end{thm}
\textbf{Proof:}\\
		Define the error term as $E(t)=x(t)-y(t)$. The synchronization is achieved if $\lim_{t \to \infty} E(t)=0$.\\
		Subtracting (\ref{4}) from (\ref{1}), we get
		\begin{equation}
			\begin{split}
				E(t) & = E(0)-E(j-1)+\frac{1}{\Gamma(\alpha)} \sum_{j=1}^t \frac{\Gamma(t-j+\alpha)}{\Gamma(t-j+1 )} \times  \\ &  \times \left[\frac{\mu(x(j-1)-y(j-1))-\mu(x^2(j-1)-y^2(j-1))}{(1+r \mu x(j-1) (1-x(j-1)))(1+r \mu y(j-1) (1-y(j-1)))} \right. \\
				& \left. -H(x(j-1),y(j-1))\right]. \label{5}
			\end{split}
		\end{equation}
		(a) Now, substituting $H=H_1$ in Eq. (\ref{5}) and taking $p=\mu$, we get
		\begin{equation*}
			\begin{split}
				E(t) &= E(0)+\frac{1}{\Gamma(\alpha)} \sum_{j=1}^t \frac{\Gamma(t-j+\alpha)}{\Gamma(t-j+1 )} [-k E(j-1)-E(j-1)].
			\end{split}
		\end{equation*}   
		From Theorem \ref{thm1},\\
		The above linear system is stable if $-1<k<-1+2^{\alpha}$.   
		
		(b) Similarly, substituting $H=H_2$ in Eq. (\ref{5}) and taking $p=\mu$, we get
		\begin{equation}
			\begin{split}
				E(t) &= E(0)+\frac{1}{\Gamma(\alpha)} \sum_{j=1}^t \frac{\Gamma(t-j+\alpha)}{\Gamma(t-j+1 )} [-k E(j-1)-E(j-1)\\
				& + \frac{\mu E^2(j-1)}{(1+r \mu x(j-1) (1-x(j-1)))(1+r \mu y(j-1) (1-y(j-1)))}]. \label{6}
			\end{split}
		\end{equation} 
		Linearizing (\ref{6}) near origin and using Theorem \ref{thm1}, we get the condition for the synchronization as $-1<k<-1+2^{\alpha}$.
		This proves the result. \\
	
	Consider the special case $r=0$. The master system (\ref{1}) becomes 
	\begin{equation}
		x(t)=x(0)+\frac{1}{\Gamma(\alpha)} \sum_{j=1}^t \frac{\Gamma(t-j+\alpha)}{\Gamma(t-j+1 )} \left[\mu x(j-1) (1-x(j-1))\right] \label{7}
	\end{equation}  
	and the slave system is
	\begin{equation}
		y(t)=y(0)+\frac{1}{\Gamma(\alpha)} \sum_{j=1}^t \frac{\Gamma(t-j+\alpha)}{\Gamma(t-j+1 )} [\mu y(j-1) (1-y(j-1))+H(x(j-1),y(j-1))].\label{8}
	\end{equation}  
	\begin{thm}
		The master system (\ref{7}) and  the slave system (\ref{8}) are synchronized if we have the following controllers :\\
		(a) $H_3= k(x-y)-p(x^2-y^2)$,\\
		(b) $H_4 = k(x-y)+2p(x y-x^2)$,\\ with the conditions $p=\mu$ and $\mu -1<k< \mu-1+2^{\alpha}$.
	\end{thm}  
	
	Proof is similar to Theorem \ref{thm2}.

	Note: 1) In \cite{wu2014chaos}, Wu and Baleanu considered systems (\ref{7}) and (\ref{8}), and studied the synchronization numerically. Here, we provided the theoretical bound on the parameter $k$. Note that the systems considered by Wu and Baleanu are slightly different. We consider the system $\Delta^{\alpha} x(t)=f(x(t+\alpha-1))-x(t+\alpha-1),$ where as the system in \cite{wu2014chaos} is $\Delta^{\alpha} x(t)=f(x(t+\alpha-1))$. Therefore, the condition for synchronization of the system in \cite{wu2014chaos} is $\mu < k< \mu + 2^{\alpha}$, the shifted interval.\\
	e.g. In \cite{wu2014chaos}, Wu and Baleanu considered as $\alpha=0.8$, $\mu=p=2.5$, $k=4$, $x_0=0.2$, $y_0=0.3$ and shown that the system is getting synchronized (cf. fig. 2 in \cite{wu2014chaos}). Using the above theory, the system can be synchronized for $2.5 < k< 4.2311$.\\
	2) Synchronization can also be achieved without these conditions. e.g., Wu and Baleanu \cite{wu2014chaos} considered $\alpha=1$, $\mu=3$, $p=2.9$, $k=3$, $x_0=0.2$, $y_0=0.3$ and demonstrated that the system is achieving synchronization (cf. fig. 4 in \cite{wu2014chaos}).\\
	i.e., The above conditions are sufficient but not necessary for the synchronization.

\begin{figure}[H]
		\centering
		\includegraphics[width=0.8\textwidth]{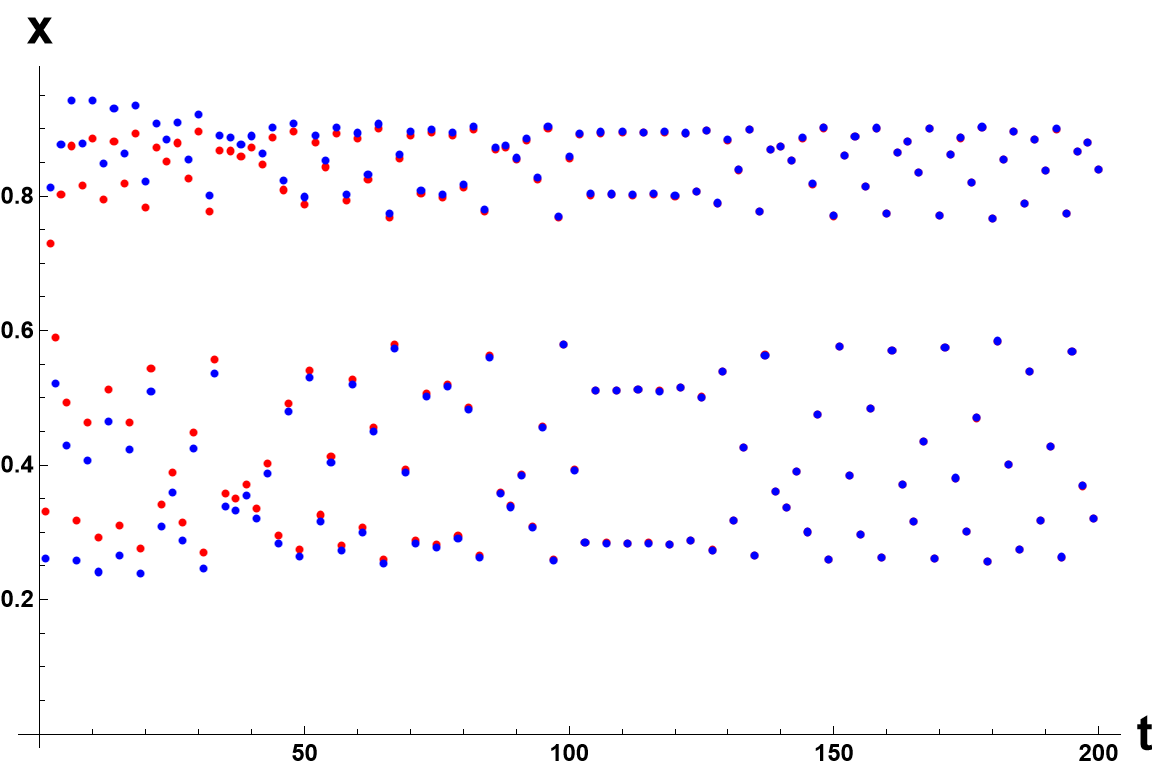}
		\caption{Synchronization of chaos for $\alpha=0.8$, $r=0.1$, $\mu=p=3.8$, and $k=0.7$. }
		\label{fig:14}
	\end{figure}

	\begin{ex}
		Consider $\alpha=0.8, \mu=3.8, r=0.1, p=3.8$ and $H=H_1$ in Theorem \ref{thm2}. In this case, $k\in(-1,0.741101)$ for synchronization. e.g. if we take $k=0.7$ and initial conditions as $x(0)=0.1$, $y(0)=0.2$ then the system (\ref{1})gets synchronized.\\
		In Fig. (\ref{fig:14}), red color represents master system and blue color represents slave system. The system (\ref{1}) and the system (\ref{4}) are synchronized after 100 iterations.
	\end{ex}

\begin{figure}[H]
		\centering
		\includegraphics[width=0.8\textwidth]{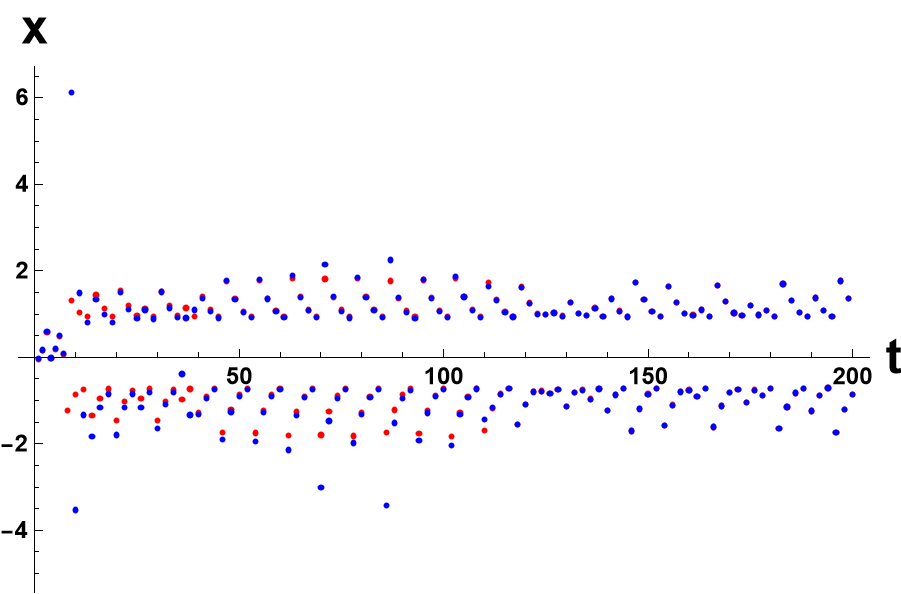}
		\caption{Synchronization of chaos for $\alpha=0.2$, $r=3.5$, $\mu=p=-3.8$, and $k=-0.5$. }
		\label{fig:15}
	\end{figure}

	\begin{ex}
		Let $\alpha=0.2, \mu=-3.8, r=3.5, p=-3.8$ and $H=H_2$ in Theorem \ref{thm2}. We need $k\in(-1,0.148698)$ to achieve the synchronization. The fig. (\ref{fig:15}) shows that systems (\ref{1}) and (\ref{4}) are in synchronized with each other after 150 iterations for $k=-0.5$ and initial conditions as $x(0)=0.01$, $y(0)=0.04$ .
	\end{ex}

	\section{Multistability} \label{mult.}
	In this section, we study the multistability in the system (\ref{1}). For $\alpha=1$, the system (\ref{1}) reduces to
	\begin{equation}
		x(t+1)=f(x(t)), \label{9}
	\end{equation}
	where $t=0,1,2,3,...$. This is a classical integer order case.    
	\par The equilibrium points of the system (\ref{9}) are the same as those of system (\ref{1}). However, the stability condition reduces to $-1<f'(x^{*})<1$ for the stability of equilibrium point $x^{*}$. Therefore, the stability region of $x_{1}^{*}$ is $-1<\mu<1$. Similarly, we can sketch the stable regions for $x_{2}^{*}$ and $x_{3}^{*}$ as in Fig. (\ref{fig2}). 
	\par
\begin{figure}[H]
		\centering
		\includegraphics[width=0.8\textwidth]{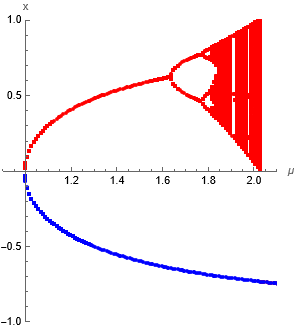}
		\caption{Multistability in the system (\ref{9}) for $r=-0.97$ and initial conditions as $x(0)=0.9$ (red color) and $x(0)=1.1$ (blue color), respectively. }
		\label{fig:16}
	\end{figure}

 The Fig. \ref{fig:16}  shows the bifurcation diagrams of the map (\ref{9}) for $r=-0.97$, $\mu \in [0,2]$ and two nearby initial conditions viz. $x(0)=0.9$ (red color) and $x(0)=1.1$ (blue color), respectively. We observed that $x(0)=0.9$ leads to period doubling and chaos whereas $x(0)=1.1$ leads to the trajectory converging to equilibrium point  $x_{3}^{*}$. 
	
\begin{figure}[H]
		\centering
		\includegraphics[width=0.8\textwidth]{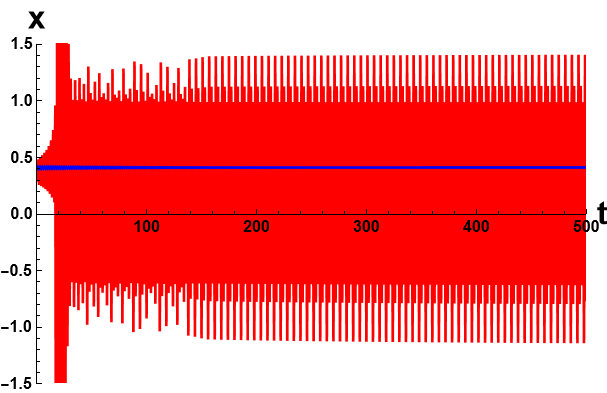}
		\caption{Multistability in the system (\ref{9}) for $\alpha=0.2$, $r=3.5$, $\mu=-3.8$, and initial conditions as $x(0)=0.3$ (red color) and  $x(0)=0.4$ (blue color), respectively. }
		\label{fig:17}
	\end{figure}

\begin{figure}[H]
		\centering
		\includegraphics[width=0.8\textwidth]{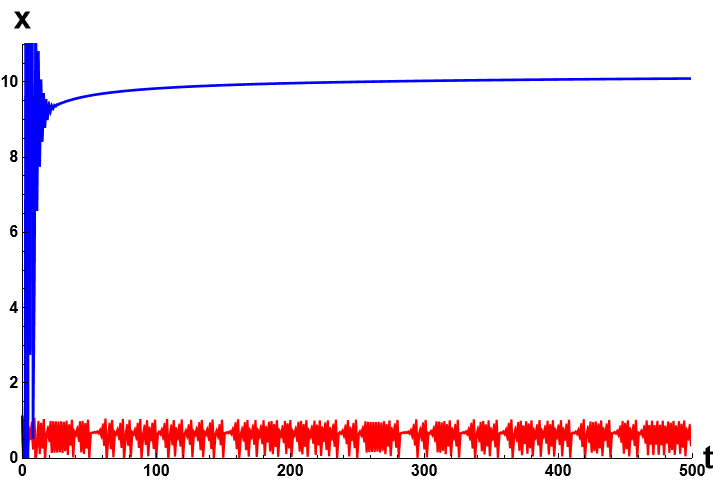}
		\caption{Multistability in the system (\ref{9}) for $\alpha=0.5$, $r=0.1$, $\mu=3.5$, and initial conditions as $x(0)=1$ (red color) and  $x(0)=1.1$ (blue color), respectively. }
		\label{fig:18}
			\end{figure}

\begin{figure}[H]
			\centering
			\includegraphics[width=0.8\textwidth]{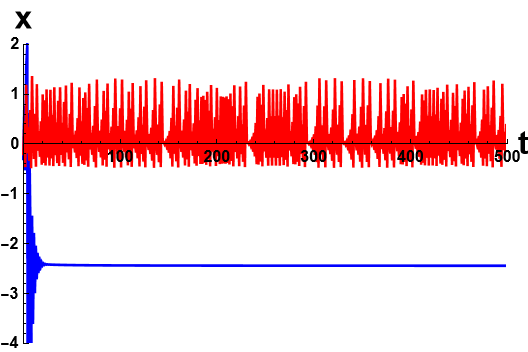}
			\caption{Multistability in the system (\ref{9}) for $\alpha=0.8$, $r=-0.5$, $\mu=-1.3$, and initial conditions as $x(0)=-0.2$ (red color) and  $x(0)=-0.3$ (blue color), respectively.}
			\label{fig:19}
		\end{figure}

	We also observed similar behavior for smaller values of $\alpha$. The Fig. \ref{fig:17} shows different behaviors of system (\ref{1}) with $r=3.5$, $\mu=-3.8$, $\alpha=0.2$ and initial conditions as $x(0)=0.3$ (red color) and  $x(0)=0.4$ (blue color), respectively. Similarly, the different behaviors  for nearby initial conditions are shown in Fig. \ref{fig:18} and Fig. \ref{fig:19} for $\alpha=0.5$ and $\alpha=0.8$, respectively. 
	

	\section{Conclusion} \label{con.}
	We generalized the classical logistic map to include the fractional order and one more parameter. The system has three equilibrium points viz. $x_1^*=0$, $x_2^*=\frac{\mu + r \mu -
		\sqrt{\mu} \sqrt{4 r + \mu - 2 r \mu+ r^2 \mu}}{2 r \mu}$ and  $x_3^*=\frac{\mu + r \mu +
		\sqrt{\mu} \sqrt{4 r + \mu - 2 r \mu+ r^2 \mu}}{2 r \mu}$. The stability of equilibrium $x_1^*$ is independent of the parameter $r$ and is stable for $1-2^{\alpha} < \mu <1 $. We provided the stability regions of $x_2^*$ and $x_3^*$ in the $\mu r$ plane. The map shows the period-doubling route to chaos. We sketched the bifurcation diagrams for various parameter values. The chaos in the proposed system is controlled by using the delayed feedback $bx(t-\tau)$. The region to achieve the control is given in the $ba$ plane using theoretical results. We also provide the results to achieve the master-slave type synchronization in the generalized logistic map. The domain of initial conditions allowed for the classical logistic map was only the unit interval. On the other hand, the proposed generalization expands the domain to the entire real line. However, the dynamic behavior can differ for the different initial conditions.

	\section*{Acknowledgments}
	P. M. Gade thanks DST-SERB for financial assistance (Ref. CRG/2020/003993). Ch. Janardhan thanks University Grants Commission, New Delhi, India for financial support (No. F.14-34/2011(CPP-II)).

	%
	

	
	\bibliographystyle{plain}      
	\bibliography{ref.bib}

\end{document}